\documentclass[graybox]{svmult}

\usepackage{type1cm}        
\usepackage{makeidx}         
\usepackage{graphicx}        
\usepackage{multicol}        
\usepackage[bottom]{footmisc}

\usepackage{newtxtext}       %
\usepackage[varvw]{newtxmath}       

\usepackage{cleveref}

\definecolor{bostonuniversityred}{rgb}{0.8, 0.0, 0.0}

\usepackage{tikz}
\usetikzlibrary{matrix}
\usetikzlibrary{shapes,patterns,snakes,arrows}
\usetikzlibrary{positioning,fit,calc}
\usetikzlibrary{graphs}

\tikzset{%
  every neuron/.style={
    circle,
    draw,
    minimum size=0.6cm
  },
  every input neuron/.style={
    circle,
    draw,
    minimum size=0.6cm,
    fill=green!50
  },
  every output neuron/.style={
    circle,
    draw,
    minimum size=0.6cm,
    fill=orange!30
  },
  every hidden neuron/.style={
    circle,
    draw,
    minimum size=0.6cm,
    fill=blue!40
  },
  neuron missing/.style={
    draw=none, 
    scale=1.5,
    text height=0.3cm,
    execute at begin node=\color{black}$\vdots$
  },
}

\makeindex             

\begin{document}

\title*{Learning Adaptive Constraints in Nonlinear FETI-DP Methods}
\author{Axel Klawonn\orcidID{0000-0003-4765-7387}\\ Martin Lanser\orcidID{0000-0002-4232-9395}\\ Janine Weber\orcidID{0000-0002-6692-2230}}
\institute{Axel Klawonn, Martin Lanser, Janine Weber \at Department of Mathematics and Computer Science, University of Cologne, Weyertal 86-90,\\ 50931 K\"oln, Germany, \url{https://www.numerik.uni-koeln.de}\\
Center for Data and Simulation Science, University of Cologne, 50923 K\"oln, Germany, \url{https://www.cds.uni-koeln.de}\\  \email{\{axel.klawonn,martin.lanser,janine.weber\}@uni-koeln.de}}
%
%
\maketitle

\abstract*{}

\abstract{
While linear FETI-DP (Finite Element Tearing and Interconnecting - Dual Primal) is an efficient iterative domain decomposition solver for discretized linear PDEs (partial differential equations), nonlinear FETI-DP is its consequent extension to the nonlinear case. In both methods, the parallel efficiency of the method results from a decomposition of the computational domain into nonoverlapping subdomains and a resulting localization of the computational work. For a fast linear convergence of the linear FETI-DP method, a global coarse problem has to be considered. Adaptive coarse spaces are provably robust variants for many complicated micro-heterogeneous problems, as, for example, stationary diffusion problems with large jumps in the diffusion coefficient. Unfortunately, the set-up and exact computation of adaptive coarse spaces is known to be computationally expensive. Therefore, recently, surrogate models based on neural networks have been trained to directly predict the adaptive coarse constraints. Here, these learned constraints are implemented in nonlinear FETI-DP and it is shown numerically that they are able to improve the nonlinear as well as linear convergence speed of nonlinear FETI-DP.    
}

\section{Introduction}
\label{sec:intro}

Linear and nonlinear FETI-DP (Finite Element Tearing and Interconnecting - Dual Primal) methods proved to be efficient parallel iterative domain decomposition solvers for linear and, respectively, nonlinear partial differential equations. While in linear FETI-DP adaptive coarse spaces (second levels) provably improve the linear convergence, that is, reduce the number of necessary Krylov iterations, it has been observed empirically that the same coarse spaces can also improve the nonlinear and linear convergence of nonlinear FETI-DP methods~\cite{MR4608253}. Unfortunately, the computation of adaptive coarse spaces is quite expensive, since many eigenvalue problems have to be set up and solved in advance. In earlier work, we have suggested two different machine learning approaches to reduce the computational effort of linear adaptive FETI-DP. In the first approach, a classification model was used to reduce the number of necessary eigenvalue problems~\cite{MR4041711}, while in the second approach a regression model was suggested to directly predict the complete coarse space without solving any eigenvalue problems~\cite{KLW_2023_learnDD}. Both models were trained for the linear stationary diffusion equation with jumps and high contrasts in the diffusion coefficient. Now, we suggest to use the trained model to predict the coarse space for Nonlinear-FETI-DP-2 and, for the first time, further combine both trained models (classification and regression model) to reduce the size of the predicted coarse space. For our numerical examples, the learned coarse space is competitive with the original adaptive approach. Let us remark that we choose Nonlinear-FETI-DP-2 as one example from the large nonlinear FETI-DP family for simplicity but similar results can be obtained, for example, for Nonlinear-FETI-DP-Res~\cite{ANE:2021:Klawonn}. 

\section{Nonlinear FETI-DP}
\label{sec:feti} 
Nonlinear FETI-DP methods are nonlinear domain decomposition methods designed to solve the discrete nonlinear equation
\begin{equation}
A(u)=0,	
\label{eq:basic}
\end{equation}
which is usually obtained by a discretization of a nonlinear partial differential equation on a domain $\Omega \in \mathbb{R}^d,\; d=2,3,$ in a finite element space $V$.  
In this section, we provide a brief description of Nonlinear-FETI-DP-2 (NL-FETI-DP-2), first introduced in~\cite{NFD:2014:Klawonn}. Let us remark that NL-FETI-DP-2 is a special case of a more general nonlinear FETI-DP framework introduced in \cite{NFD:2017:Klawonn} and extended in \cite{ANE:2021:Klawonn}. All types of learned or computed coarse constraints discussed here can also be used in the general framework, but we decided to concentrate on NL-FETI-DP-2 due to page limitations. For 
more details on nonlinear FETI-DP methods, see~\cite{NFD:2014:Klawonn,TES:2015:Klawonn,NFD:2017:Klawonn}. 
In FETI-DP, we consider nonoverlapping subdomains $\Omega_i,\,i=1,...,N$ with $\bigcup_{i=1}^{N} \Omega_i =\Omega$ and corresponding local finite element spaces $W_i,\,i=1,...,N$. We now introduce local nonlinear problems 
$$
K_i(u_i) -f_i =0,\ i=1,...,N,
$$ 
on the subdomains, which are obtained by a local finite element discretization, assuming zero Neumann 
type boundary conditions on the interface 
\begin{equation} \label{eq:interface}
\Gamma = \bigcup_{i \neq j} \left( \partial \Omega_i \cap \partial \Omega_j \right) \setminus \partial \Omega_D.	
\end{equation}
By $\partial \Omega_D \subseteq \partial \Omega$ we denote the part of the boundary where Dirichlet type boundary conditions are imposed. We can always compute the original problem~\cref{eq:basic} by
$$
A(u) = R^TK(Ru) -R^Tf,
$$
where $W= W_1 \times \cdots \times W_N$, $R_i:V \rightarrow W_i$, $R=\left(R_1^T,...,R_N^T\right)^T$, $f=\left(f_1^T,...,f_N^T\right)^T$, and $K(Ru) = \left(K_1(R_1u)^T,...,K_N(R_Nu)^T\right)^T$. 
For the FETI-DP method, all interface variables are now partitioned into dual variables (index set $\Delta$) and primal variables (index set $\Pi$). The primal variables or primal constraints can be subdomain vertices as well as (weighted) edge or face averages. In this article, we will consider adaptive and learned constraints in 2D, which can be classified as weighted edge constraints. The primal constraints can be interpreted as the coarse constraints of FETI-DP and build a global coarse problem or, in other words, second level. We now define the space $\widetilde{W}$ of functions, which are only enforced to be continuous in all primal variables. With the restriction $\check{R}: W \rightarrow \widetilde{W}$, $\tilde{u} \in \widetilde{W}$, we can define 
$$
\widetilde{K}(\tilde{u}) := \check{R}^T K(\check{R}\tilde{u}), \; \tilde{f} := \check{R}^T f. 
$$
Additionally, we also have to control the jump in the dual variables and therefore introduce the jump matrix $B: W \rightarrow {\rm range}(B)$, which computes the jump across the interface of functions from $W$ (see~\cite{HSP:2010:Klawonn,toselli:2004:ddm} for a detailed definition). Enforcing the constraint $B\check{R}\tilde{u}=0$ ensures that $\tilde{u}$ is continuous on the complete interface. Therefore, we introduce Lagrange multipliers $\lambda$, and define the nonlinear saddle point problem
\begin{equation}
\left(\begin{array}{l}
	\widetilde{K}(\tilde{u}) + \check{R}^T B^T \lambda - \tilde{f}\\
	B\check{R}\tilde{u}
\end{array}\right) = \left(\begin{array}{l}
	0\\0
\end{array}\right).
\label{eq:sps}
\end{equation}
All nonlinear FETI-DP methods are based on solving this saddle point problem, that is, the search for the unique solution $(\tilde{u}^*, \lambda^*)$. In general, an arbitrary set of variables can be eliminated nonlinearly and the resulting nonlinear Schur complement is solved with Newton's method. Here, we focus on NL-FETI-DP-2, where the complete first line in \cref{eq:sps} is eliminated nonlinearly. That is, we reformulate the first equation to
\begin{equation}
\tilde{u} = \widetilde{K}^{-1}(\check{R}^T B^T \lambda - \tilde{f})
\label{eq:fl}
\end{equation}
and afterwards solve the nonlinear Schur complement 
$$
B\check{R} \widetilde{K}^{-1}(\check{R}^T B^T \lambda - \tilde{f}) =0
$$ with Newton's method for $\lambda^*$. In each Newton step $\lambda^{(k)}$, \cref{eq:fl} has to be solved with $\lambda := \lambda^{(k)}$, which is also done using Newton's method. Therefore, we have two nested Newton iterations in our method. In the inner loop, an equation system with the Jacobian matrix $D\widetilde{K}(\cdot)$ has to be solved. Due to the special structure of the matrix, this can be done efficiently and in parallel by local sparse direct solvers on each subdomain which eliminate the dual and interior variables and an additional global sparse direct solver for the resulting coarse problem in the primal variables. We refer to~\cite{TES:2015:Klawonn} for details. In the outer loop, the Jacobian matrix writes
$$
F(\cdot) := B\check{R} D\widetilde{K}(\cdot) \check{R}^T B^T,
$$ 
which has the same structure as the usual linear FETI-DP matrix $F$; see~\cite{NFD:2014:Klawonn} for details. The linearized system can thus be solved iteratively with a Krylov subspace method and using the usual preconditioner
\begin{equation}
M^{-1}_D = \sum_{i=1}^N B_D^{(i)} S^{(i)} B_D^{(i)T}.
\label{eq:dir_prec}
\end{equation}  
Here, $S^{(i)}$ is the Schur complement of $DK_i(\cdot)$ with respect to the interface, where $DK_i(\cdot)$ is the Jacobian matrix of $K_i(\cdot)$.
 Finally, $B_D^{(i)}$ is the scaled local jump matrix $B^{(i)}$, with  $B = (B^{(1)},...,B^{(N)})$ and $B_{D} = (B_D^{(1)},...,B_D^{(N)})$. For a description of different scalings, we refer to~\cite{DPF:2006:Klawonn,ACO:2016:Klawonn}. Throughout this paper, we exclusively use $\rho$-scaling. 
 
 \section{Adaptive and learned FETI-DP coarse spaces}
\label{sec:adaptive-feti}

To obtain a robust FETI-DP method, appropriate primal constraints have to be chosen. For many linear problems, adaptive coarse spaces are provably robust, 
but require the solution of one generalized eigenvalue problem for each edge of the domain decomposition. In this article, we use the adaptive coarse space introduced in~\cite{sousedik1}. We will not introduce the eigenvalue problem here due to page limitations. Let us only remark that for, e.g., linear diffusion or elasticity problems, the condition number $\kappa$ of the preconditioned system $M_D^{-1} F$ is then bounded by
$$\kappa \leq N_{E}^{2} TOL,$$ 
where $TOL$ is a tolerance given by the user and $N_{E}$ a small geometric constant. Depending on the value of $TOL$, the first $k$ eigenvectors belonging to the eigenvalues which are larger than $TOL$ on an edge are used to create $k$ constraints. Let us remark that on many edges $k$ is equal to 0 and, hence, 0 constraints are necessary for this edge. For the eigenvalue problem on an edge, the local stiffness matrices of the two adjacent subdomains $\Omega_i$ and $\Omega_j$ are needed. In the nonlinear case, we always use the Jacobian matrices $DK_i(\cdot)$ and $DK_j(\cdot)$ instead, evaluated in the initial value of Newton's method. The constraints are always only computed once and then fixed for the complete nonlinear solve. 

To avoid the solution of many eigenvalue problems in the linear case, we have set up two different approaches based on neural networks. The first one is a classification approach~\cite{MR4041711} which predicts how many adaptive constraints are needed on an edge. It classifies the edges into the classes $k=0$, $k=1$, and $k>1$ and only for the latter class eigenvalue problems have to be solved. The second approach is a regression network which directly predicts the first three constraints on each edge~\cite{KLW_2023_learnDD}. Note that actually six different regression networks are trained to predict the first three constraints, that is, two for each of the three constraints depending on the connection of the edge with the Dirichlet boundary; see~\cite{KLW_2023_learnDD} for more details. In general, using this approach, all eigenvalue problems can be avoided.

\begin{figure}[t]
\centering
\scalebox{0.68}{
\begin{tikzpicture}[line width=1pt]

\path[draw] (-3,0) -- (3,0) -- (3,3) -- (-3,3) -- cycle;

\path[draw,fill=bostonuniversityred!30] (-0.6,0.6) -- (1.8,0.6) -- (1.8,1.2) -- (-0.6,1.2) -- cycle;
\path[draw,fill=bostonuniversityred!30] (-2.4,1.2) -- (-1.8,1.2) -- (-1.8,1.8) -- (0.0,1.8) -- (0.0,2.4) -- (-2.4,2.4) -- cycle;
\path[draw,fill=bostonuniversityred!30] (1.2,1.8) -- (2.4,1.8) -- (2.4,3.0) -- (1.2,3.0) -- cycle;

\foreach \y in {0.0,0.6,...,3.0} {
	\foreach \x in {-3,-2.4,...,3} {
		\path [draw,line width=0.5pt] (\x,\y) -- (\x+0.6,\y+0.6) -- (\x,\y+0.6) -- cycle;
	}
}

\path[draw,line width=3.0pt,blue] (0,0) -- (0,3);

\foreach \y in {0,...,12}{
	\foreach \x in {1,...,12}{
		\draw[fill=bostonuniversityred] (0.25*\x,0.25*\y) circle (0.04);
	}
	\foreach \x in {-12,...,-1}{
		\draw[fill=bostonuniversityred] (0.25*\x,0.25*\y) circle (0.04);
	}
}

\draw (-1.5,2.5) node[text=black,fill=white] {$\Omega_i$};
\draw (1.5,2.5) node[text=black,fill=white] {$\Omega_j$};
\draw (0,-0.5) node[text=blue,fill=black!0] {$\mathcal{E}_{ij}$};

\def\xoffset{4.2}
\def\yoffset{0.35}

\foreach \m/\l [count=\y] in {1,2,missing,3,missing,4}
  \node [every neuron/.try, neuron \m/.try] (input-\m) at (\xoffset+0,\yoffset+3.75-0.75*\y) {};

\foreach \m [count=\y] in {1,2,missing,3}
  \node [every neuron/.try, neuron \m/.try ] (hidden1-\m) at (\xoffset+1.6,\yoffset+3.5-\y*1.0) {};
  
\foreach \m [count=\y] in {1,2,missing,3}
  \node [every neuron/.try, neuron \m/.try ] (hidden2-\m) at (\xoffset+3.2,\yoffset+3.5-\y*1.0) {};  

\foreach \m [count=\y] in {1,missing,2,missing,3}
  \node [every neuron/.try, neuron \m/.try ] (output-\m) at (\xoffset+4.8,\yoffset+3.4-0.75*\y) {};

\foreach \i in {1,...,4}
  \foreach \j in {1,...,3}
    \draw [->] (input-\i) -- (hidden1-\j);
    
\foreach \i in {1,...,3}
  \foreach \j in {1,...,3}
    \draw [->] (hidden1-\i) -- (hidden2-\j);

\foreach \i in {1,...,3}
  \foreach \j in {1,...,3}
    \draw [->] (hidden2-\i) -- (output-\j);

\foreach \l [count=\x from 0] in {Input\\layer, Hidden\\layers, Output\\layer, Learned\\constraint}
  \node [align=center, above] at (\xoffset+\x*2.4,4.0) {\l};

\draw[->,bostonuniversityred,smooth] (0.25,3.0) .. controls (2.0,3.5) .. (input-1) 
                     node [fill=white, above=0.1cm, pos=0.9] {$\rho(x_1)$};
\draw[->,bostonuniversityred,smooth] (0.25,2.75) .. controls (1.75,1.6) .. (input-2) 
                     node [fill=white, above=0.3cm, pos=0.89] {$\rho(x_2)$};
\draw[->,bostonuniversityred,smooth] (-2.0,2.0) .. controls (0.5,1.0) .. (input-3) 
                     node [fill=white, below=0.1cm, pos=0.94] {$\rho(x_k)$};
\draw[->,bostonuniversityred,smooth] (3.0,0.0) .. controls (3.3,-0.2) .. (input-4) 
                     node [fill=white, below=0.4cm, pos=0.1] {$\rho(x_M)$};

\node (constr_pic) at (11.8,1.65) {\includegraphics[width=.3\textwidth]{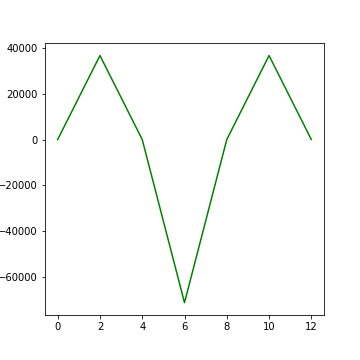}};

\draw [->,bostonuniversityred,smooth]  (output-1) .. controls (10,2.3) .. (10.4,2.1);
\draw [->,bostonuniversityred,smooth]  (output-2) .. controls (10,1.3) .. (11.7,0.5);
\draw [->,bostonuniversityred,smooth]  (output-3) .. controls (11.8,-0.9) .. (13.2,2.1);

\end{tikzpicture}
}
\caption{Visualization of our regression neural network model. As input data for the neural network,
we use samples of the coefficent function for the two neighboring subdomains of an edge (left).
Here, dark red corresponds to a high coefficient and white corresponds to a low coefficient. The
output of the network is a discretized egde constraint (right). Figure taken from~\cite[Fig. 1]{KLW_2023_learnDD}.}
\label{fig:fnn}
\end{figure}
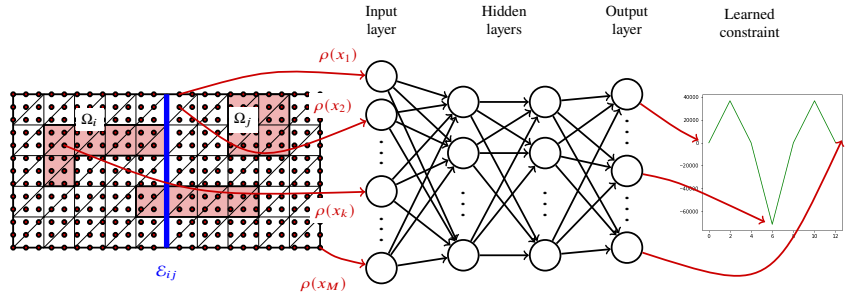

Both networks are trained for linear diffusion problems; see~\cite{MR4041711,KLW_2023_learnDD} for details. 
In particular, both approaches, that is, the classification and regression neural networks, are trained using local function evaluations of the coefficient or material distribution of the underlying test problem within the two adjacent subdomains of an edge. Hence, the training is completely localized in the sense that, within one evaluation of the neural networks, only a specific edge as well as its directly neighboring subdomains are considered. Moreover, the training of the neural networks is independent of the finite element discretization as long as the discretization is fine enough to capture all heterogeneities within the coefficient function. 
As output data of the classification neural network, we use the information from the eigenvalue problem whether for the respective edge, $k=0$, $k=1$, or $k>1$ adaptive constraints are necessary to obtain a robust coarse space.
For the regression network, we use discrete approximations of the first three adaptive edge constraints resulting from the eigenvalue problem defined in~\cite{sousedik1} as output data. Let us note that due to the fixed number of output nodes of the trained neural networks, in principle, the output space of the regression neural network corresponds to a fixed edge length. Hence, in order to allow the application of our approach for different mesh sizes, we use a linear interpolation of the discretized edge constraints using the finite element mesh points as the interpolation points and the finite element basis functions as interpolation basis.
A schematic representation of the trained regression neural networks can be seen in~\cref{fig:fnn}.

In the following, we simply use both networks trained for linear FETI-DP to construct learned coarse spaces for NL-FETI-DP-2 for the solution of nonlinear $p$-Laplace equations.  
In particular, for the first time, we combine both networks such that we first evaluate the classification neural network to obtain an a priori estimation of the number of necessary coarse constraints and, subsequently, implement the respective number of learned constraints as predicted by the regression neural network. 

\section{Test Problem and Numerical Results}
\label{sec:results}

\begin{figure}
\centering
\includegraphics[width=0.45\textwidth]{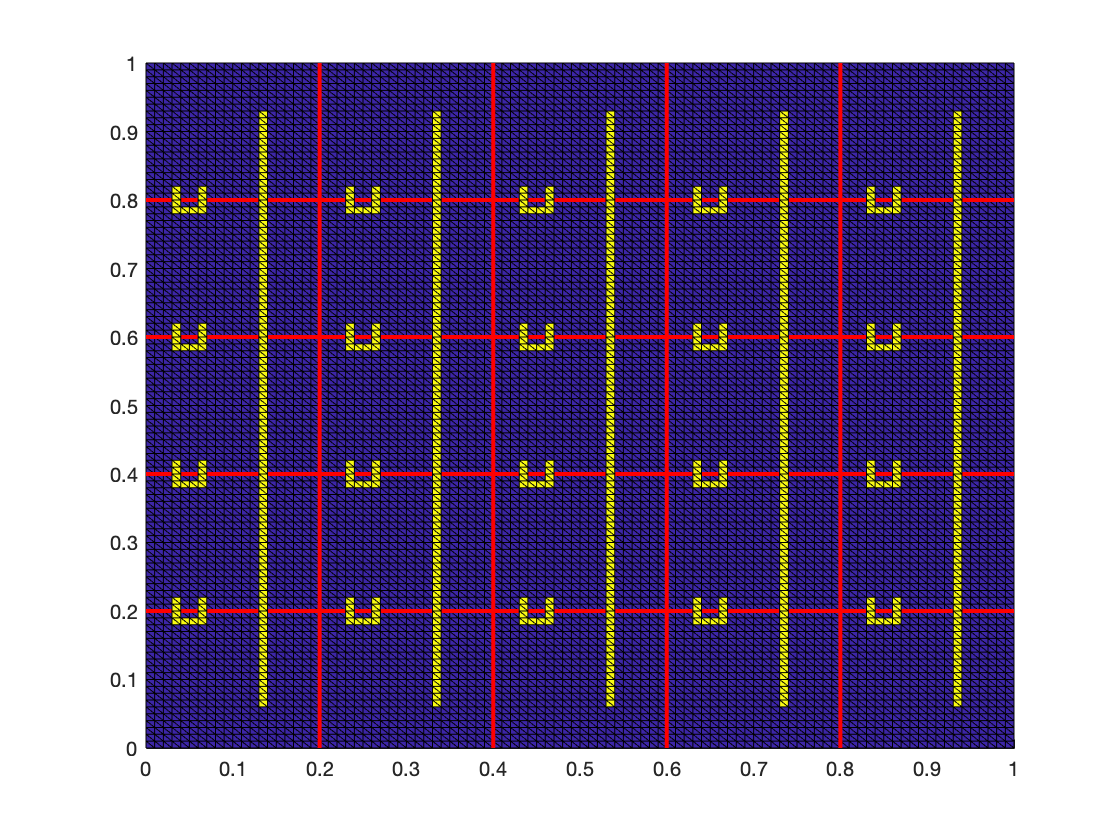}
\includegraphics[width=0.45\textwidth]{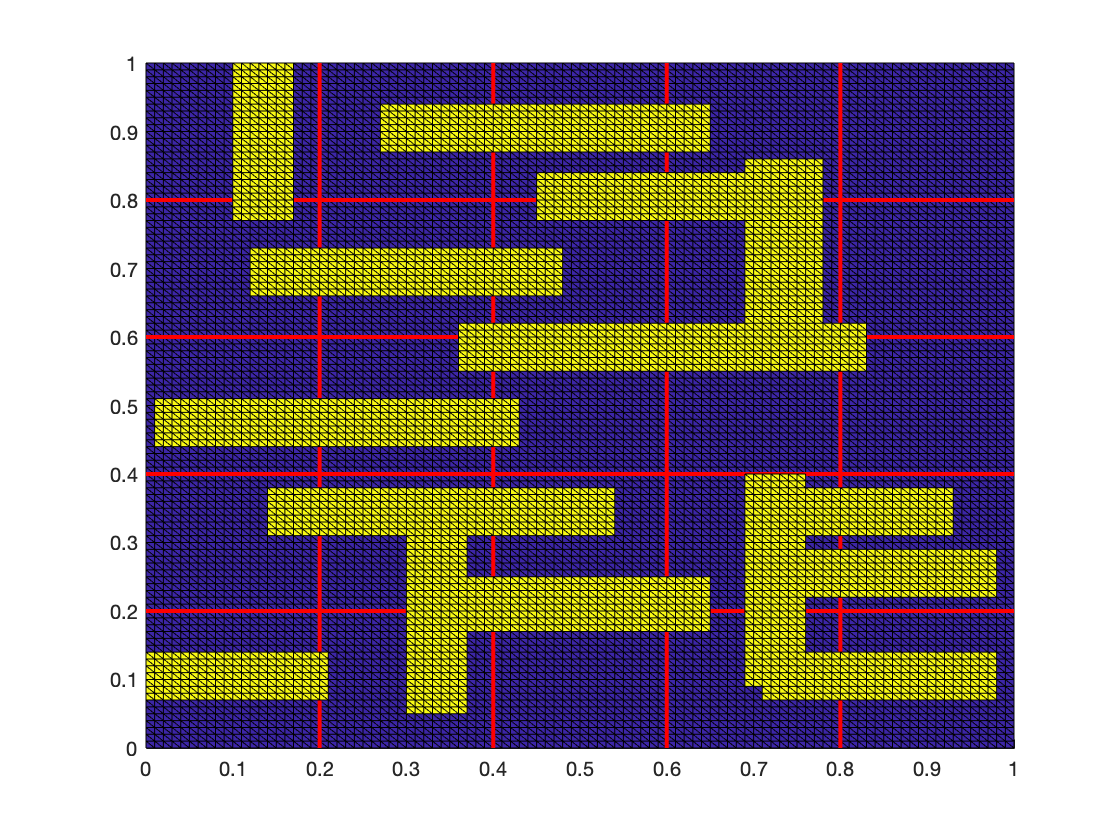}
\caption{Two different coefficient distributions $\alpha$, where $\alpha$ is equal to one in the blue finite elements and $1e6$ in the yellow ones.}
\label{fig:alpha}	
\end{figure}

We consider the model problem
\begin{equation}
\begin{array}{rccl}
- \alpha \Delta_{p} u  &=&  1 &\qquad {\rm in }\; \Omega, \\ 
 u          &=&  0 &\qquad {\rm on } \;\partial \Omega_D,
\end{array}
\label{lanser_MS06:klawonna_plenary_3:laplace}
\end{equation}
with the scaled $p$-Laplace operator
$\alpha \Delta_{p} u := {\rm div} (\alpha |\nabla u |^{p-2} \nabla u)$.
Within this article, we use $p=4$ and a coefficient function
$\alpha: \Omega \rightarrow \mathbb{R}$ with jumps. 
We always use the unit square $\Omega = [0,1] \times [0,1]$ as the computational domain, a discretization with piecewise linear finite elements, and a structured domain decomposition into square subdomains. Let us note that $\partial \Omega_D$ is chosen to be the left boundary ($x=0$) of the domain and that we enforce natural Neumann-type boundary conditions on the remaining boundary. In our numerical experiments, we consider two different coefficient functions $\alpha$ and a decomposition into $5 \times 5$ subdomains; see \cref{fig:alpha}. As initial value for Newton's method, we always use $u^{(0)}(x,y)=x$. 

We combine classical Newton-Krylov-FETI-DP (NK-FETI-DP), which is applying Newton's method to solve the global problem in~\cref{eq:basic} and then using linear FETI-DP for all linear solves, and Nonlinear-FETI-DP-2 (NL-FETI-DP-2) with four different coarse spaces, that is, using only vertices, adaptive constraints computed with $TOL=100$, two learned constraints on each edge, and using the classification network to decide how many learned constraints are used on each edge. In the latter approach, we either use no edge constraint (class 0), the first learned constraint (class 1), or the first and second learned constraint (class 2) on each edge. 

All the results are presented in~\cref{tab:res} where several things can be observed.  
For the first example, in principle, also vertices can be used in NL-FETI-DP-2 and result in a fast nonlinear and linear convergence. Nonetheless, the adaptive as well as the learned constraints perform better and the quality of both are nearly identical. For NK-FETI-DP, one should avoid the vertex coarse space due to a slow linear convergence. 

For the second example, we observe more differences. In the case of NK-FETI-DP, the learned coarse spaces have larger condition numbers compared with the adaptive coarse space, but the linear convergence is still convincing. The nonlinear convergence of NL-FETI-DP-2 with a vertex coarse space is really slow which is improved drastically using adaptive or learned coarse spaces. The learned coarse spaces are as good as the adaptive coarse space and clearly can compete. In general, in both examples, the combination of the classification and the regression model yields a smaller size of the coarse space without deteriorating the convergence.

\begin{figure}
\centering
\includegraphics[width=0.75\textwidth]{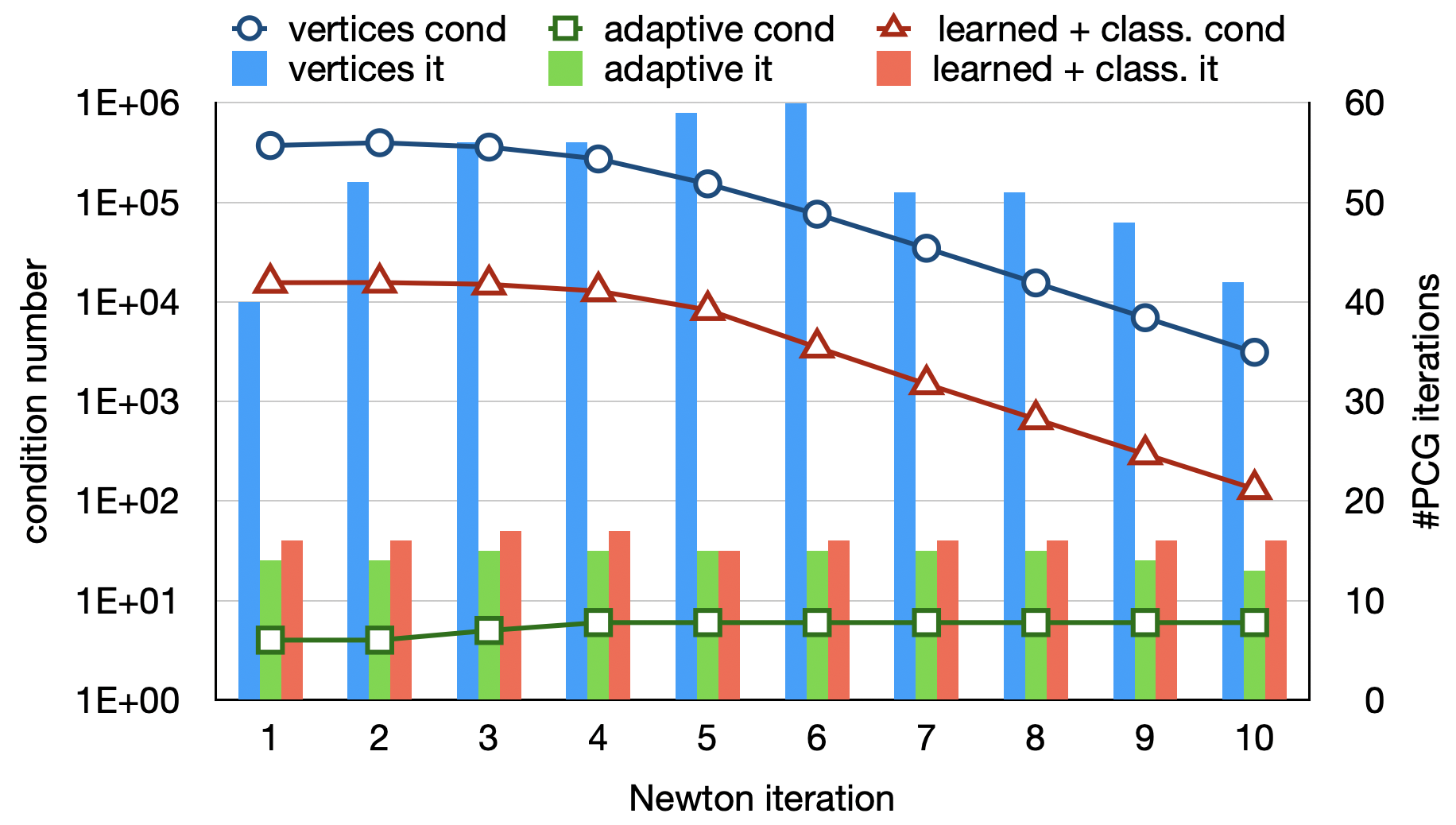}
\caption{Condition number and PCG iterations in each Newton step of NK-FETI-DP with three different coarse spaces (vertices, adaptive with $TOL=100$, and classification plus learned approach) for the example with the coefficient function shown in~\cref{fig:alpha} (right).}
\label{fig:its}	
\end{figure}

Finally, in \cref{fig:its}, we present the development of the condition number and PCG iteration counts over the Newton steps of NK-FETI-DP. While the condition number is small and stays nearly constant using the adaptive coarse space, it is much higher and tends to decrease for the other two cases. The  condition number implementing the learned coarse space is two orders of magnitude smaller than the condition number using solely vertex constraints and the PCG iteration counts are comparable with the adaptive case.

All in all, the learned constraints improve the nonlinear and linear convergence of NL-FETI-DP-2 as much as the adaptive coarse space and can be used in the nonlinear case without any modifications.

\begin{table}[t]
\begin{center}
\small
\caption{Results for both coefficient distributions; four different coarse spaces used in NK-FETI-DP and NL-FETI-DP-2; {\bf outer/inner/PCG it.} stand for the Newton iterations in the outer loop, the inner Newton iterations summed over the outer ones, and the PCG iterations summed over the outer iterations; $|\Pi|$ gives the size of the respective coarse space; {\bf min./max. cond.} gives the minimal and maximal estimated condition number in all outer Newton steps.}
\begin{tabular}{|ll|r|rrr|rr|}
\hline 
\multicolumn{8}{|c|}{\bf Channels and Us; see~\cref{fig:alpha} (left)} \\
\multicolumn{8}{|c|}{ $p = 4$; $H/h=20$; 25 subdomains}\\\hline
\bf  & \bf coarse & \bf & \bf outer & \bf inner & \bf PCG & \bf min.  & \bf max. \\
\bf method & \bf space & \bf $|\Pi|$ & \bf it. & \bf it. & \bf it. & \bf cond.  & \bf cond. \\\hline
NK-FETI-DP&vertices&28&10&-&781&195.6&42\,690.1\\
NK-FETI-DP&adaptive; $TOL=100$&68&10&-&99&4.9&5.5\\
NK-FETI-DP&learned&108&10&-&63&1.3&1.4\\
NK-FETI-DP&learned + class.&68&10&-&99&4.9&5.5\\\hline 
NL-FETI-DP-2&vertices&28&2&20&46&11.9&21.0\\
NL-FETI-DP-2&adaptive; $TOL=100$&68&1&11&7&1.9&1.9\\
NL-FETI-DP-2&learned&108&1&11&6&1.4&1.4\\
NL-FETI-DP-2&learned + class.&68&1&11&7&1.9&1.9\\
\hline\hline
\multicolumn{8}{|c|}{\bf Combs; see~\cref{fig:alpha} (right)} \\
\multicolumn{8}{|c|}{ $p = 4$; $H/h=20$; 25 subdomains}\\\hline
\bf  & \bf coarse & \bf & \bf outer & \bf inner & \bf PCG & \bf min.  & \bf max. \\
\bf method & \bf space & \bf $|\Pi|$ & \bf it. & \bf it. & \bf it. & \bf cond.  & \bf cond. \\\hline
NK-FETI-DP&vertices&28&10&-&516&3\,116.3&397\,020.1\\
NK-FETI-DP&adaptive; $TOL=100$&41&10&-&145&4.3&6.8\\
NK-FETI-DP&learned&108&10&-&124&122.0&15\,083.1\\
NK-FETI-DP&learned + class.&41&10&-&161&132.9&15\,642.0\\\hline 
NL-FETI-DP-2&vertices&28&9&53&257&16.8&6\,746.7\\
NL-FETI-DP-2&adaptive; $TOL=100$&41&1&20&17&14.0&14.0\\
NL-FETI-DP-2&learned&108&2&22&27&13.0&20.0\\
NL-FETI-DP-2&learned + class.&41&1&20&18&14.0&14.0\\
\hline
\end{tabular}
\label{tab:res}
\end{center}
\end{table}

\bibliographystyle{abbrv}
\bibliography{NL-FETI-DP}{}

\begin{thebibliography}{10}

\bibitem{MR4608253}
A.~Heinlein, A.~Klawonn, and M.~Lanser.
\newblock Adaptive nonlinear domain decomposition methods with an application
  to the {$p$}-{L}aplacian.
\newblock {\em SIAM J. Sci. Comput.}, 45(3):S152--S172, 2023.

\bibitem{MR4041711}
A.~Heinlein, A.~Klawonn, M.~Lanser, and J.~Weber.
\newblock Machine learning in adaptive domain decomposition
  methods---predicting the geometric location of constraints.
\newblock {\em SIAM J. Sci. Comput.}, 41(6):A3887--A3912, 2019.

\bibitem{NFD:2014:Klawonn}
A.~Klawonn, M.~Lanser, and O.~Rheinbach.
\newblock Nonlinear {FETI}-{DP} and {BDDC} methods.
\newblock {\em SIAM J. Sci. Comput.}, 36(2):A737--A765, 2014.

\bibitem{TES:2015:Klawonn}
A.~Klawonn, M.~Lanser, and O.~Rheinbach.
\newblock Toward extremely scalable nonlinear domain decomposition methods for
  elliptic partial differential equations.
\newblock {\em SIAM J. Sci. Comput.}, 37(6):C667--C696, 2015.

\bibitem{NFD:2017:Klawonn}
A.~Klawonn, M.~Lanser, O.~Rheinbach, and M.~Uran.
\newblock Nonlinear {FETI}-{DP} and {BDDC} methods: a unified framework and
  parallel results.
\newblock {\em SIAM J. Sci. Comput.}, 39(6):C417--C451, 2017.

\bibitem{ANE:2021:Klawonn}
A.~Klawonn, M.~Lanser, and M.~Uran.
\newblock Adaptive nonlinear elimination in nonlinear {FETI-DP} methods.
\newblock In {\em Domain Decomposition Methods in Science and Engineering
  XXVI}, pages 337--345. Springer, 2023.

\bibitem{KLW_2023_learnDD}
A.~Klawonn, M.~Lanser, and J.~Weber.
\newblock Learning adaptive coarse basis functions of {FETI-DP}.
\newblock {\em Journal of Computational Physics}, 496:112587, 2024.

\bibitem{ACO:2016:Klawonn}
A.~Klawonn, P.~Radtke, and O.~Rheinbach.
\newblock A comparison of adaptive coarse spaces for iterative substructuring
  in two dimensions.
\newblock {\em Electron. Trans. Numer. Anal.}, 45:75--106, 2016.

\bibitem{HSP:2010:Klawonn}
A.~Klawonn and O.~Rheinbach.
\newblock Highly scalable parallel domain decomposition methods with an
  application to biomechanics.
\newblock {\em ZAMM Z. Angew. Math. Mech.}, 90(1):5--32, 2010.

\bibitem{DPF:2006:Klawonn}
A.~Klawonn and O.~B. Widlund.
\newblock Dual-primal {FETI} methods for linear elasticity.
\newblock {\em Comm. Pure Appl. Math.}, 59(11):1523--1572, 2006.

\bibitem{sousedik1}
J.~Mandel and B.~Soused{\' i}k.
\newblock Adaptive selection of face coarse degrees of freedom in the {BDDC}
  and the {FETI}-{DP} iterative substructuring methods.
\newblock {\em Comput. Methods Appl. Mech. Engrg.}, 196(8):1389--1399, 2007.

\bibitem{toselli:2004:ddm}
A.~Toselli and O.~Widlund.
\newblock {\em Domain Decomposition Methods - Algorithms and Theory}, volume~34
  of {\em Springer Series in Computational Mathematics}.
\newblock Springer, 2004.

\end{thebibliography}

\end{document}